\documentclass[11pt]{article}
\usepackage{amssymb,amsfonts,amsmath,latexsym,epsf,tikz,url}

\newtheorem{theorem}{Theorem}[section]

\newtheorem{observation}[theorem]{Observation}

\newtheorem{corollary}[theorem]{Corollary}
\newtheorem{lemma}[theorem]{Lemma}
\newtheorem{remark}[theorem]{Remark}

\newtheorem{problem}[theorem]{Problem}

\newcommand{\proof}{\noindent{\bf Proof.\ }}
\newcommand{\qed}{\hfill $\square$\medskip}

\begin{document}

\title{The chromatic distinguishing index of certain graphs}

\author{
Saeid Alikhani  $^{}$\footnote{Corresponding author}
\and
Samaneh Soltani
}

\date{\today}

\maketitle

\begin{center}
Department of Mathematics, Yazd University, 89195-741, Yazd, Iran\\
{\tt alikhani@yazd.ac.ir, s.soltani1979@gmail.com}
\end{center}

%%%%%%%%%%%%%%ABSTRACT%%%%%%%%%%%%%%%%%%%%%%%%%%%%%%%%%%%%%%%%%%%%%%%%%%%%%%%%%%%%

\begin{abstract}
   The distinguishing index of a graph $G$, denoted by $D'(G)$, is the
least number of labels in an edge labeling of $G$ not preserved by
any non-trivial automorphism.  The distinguishing chromatic index $\chi'_D (G)$ of a graph $G$ is the least number $d$ such that $G$ has a proper edge labeling with $d$ labels that is preserved only by the identity automorphism of $G$. In this paper we compute the distinguishing chromatic index for some specific graphs.   Also we study the distinguishing chromatic index of corona product and join of two graphs.  
\end{abstract}


\medskip
\noindent{\bf AMS Subj.\ Class.}: 05C25, 05C15

%%%%%%%%%%%%%%%%%%%%%%%%%%%%%%%%%%%%%%%%%%%%%%%%%%%%%%%%%%%%%%%%%%%%%%%%%%%%%%%%%
%%%%%%%%%%%%%%%%%%%%%%%%%%%%%%%%%%%%%%%%%%%%%%%%%%%%%%%%%%%%%%%%%%%%%%%%%%%%%%%%%
\section{Introduction}
%%%%%%%%%%%%%%%%%%%%%%%%%%%%%%%%%%%%%%%%%%%%%%%%%%%%%%%%%%%%%%%%%%%%%%%%%%%%%%%%%
Let $G=(V,E)$ be a simple graph and ${\rm Aut}(G)$ denote the automorphism group of $G$.  For  $v \in V$, the \textit{neighborhood} of  $v$ is the set $N_G(v) = \{u \in V(G) : uv \in   E(G)\}$. The \textit{degree} of  $v$ in a graph $G$, denoted by ${\rm deg}_G(v)$, is the number of edges of $G$ incident with $v$. In particular, ${\rm deg}_G(v)$ is the number of neighbours of $v$ in $G$.  Also, the maximum degree of $G$  is denoted by $\Delta (G)$.

A proper edge labeling $c$ of a nonempty graph $G$ (a graph with edges) is a function
$c : E(G)\rightarrow S$, where $S$ is a set of labels (colors), with the
property that $c(e)\neq c(f)$ for every two adjacent edges $e$ and $f$ of $G$. If the labels are
chosen from a set of $k$ labels, then $c$ is called a \textit{proper $k$-edge labeling} of $G$. The minimum
positive integer $k$ for which $G$ has a proper $k$-edge labeling is called the \textit{chromatic index} of
$G$ and is denoted by $\chi'(G)$.
As a result of Vizing's theorem, the chromatic index of every nonempty graph $G$
is one of two numbers, namely $\Delta(G)$ or $\Delta(G)+1$. A graph $G$ with $\chi'(G)= \Delta(G)$
is called a \textit{class one} graph while a graph $G$ with $\chi'(G)= \Delta(G)+1$ is called a \textit{class
two graph}. For instance, it is proved that, $K_n$ is a class one graph if $n$ is even and is a class two graph if $n$ is odd, and also  every regular graph of odd order is a class two graph.
The next two results is about graphs whose  are class one. 
\begin{theorem}{\rm \cite{69}}\label{thm1.5}
 Every bipartite graph is a class one graph.
\end{theorem}
\begin{corollary}{\rm \cite{43}}\label{mohemmaximum}
 If $G$ is a graph in which no two vertices of maximum degree
are adjacent, then $G$ is a class one graph.
\end{corollary}

A proper vertex labeling of a graph $G$ is a function $c : V(G)\rightarrow S$, such that $c(u)\neq c(v)$ for every pair $u$ and $v$ of adjacent vertices of $G$. If $|S|=  k$, then $c$ is called a \textit{proper $k$-vertex labeling} of $G$. The minimum positive integer $k$ for which
G has a proper $k$-vertex labeling is called the \textit{chromatic number}  of $G$ and  is denoted by $\chi(G)$.

A labeling of $G$, $\phi : V \rightarrow \{1, 2, \ldots , r\}$, is said to be \textit{$r$-distinguishing}, 
if no non-trivial  automorphism of $G$ preserves all of the vertex labels.
The point of the labels on the vertices is to destroy the symmetries of the
graph, that is, to make the automorphism group of the labeled graph trivial.
Formally, $\phi$ is $r$-distinguishing if for every non-trivial $\sigma \in {\rm Aut}(G)$, there
exists $x$ in $V$ such that $\phi(x) \neq \phi(\sigma(x))$. Authors  often refer to a
labeling as a coloring, but there is no assumption that adjacent vertices get
different colors. Of course the goal is to minimize the number of colors used.
Consequently  the \textit{distinguishing number}  of a graph $G$ is defined  by
\begin{equation*}
D(G) = min\{r \vert ~ G ~\textsl{has a labeling that is $r$-distinguishing}\}.
\end{equation*} 

This number has defined in \cite{Albert}. If a graph has no nontrivial automorphisms, its distinguishing number is  $1$. In other words, $D(G) = 1$ for the asymmetric graphs.
 The other extreme, $D(G) = \vert V(G) \vert$, occurs if and only if $G = K_n$.  Collins and Trenk \cite{EJC} defined the \textit{distinguishing chromatic number}
$\chi_D(G)$ of a graph $G$ for proper labelings, so $\chi_D(G)$ is the least number $d$ such that $G$ has a proper labeling with $d$ labels that is only preserved by the trivial automorphism. 
Similar to this definitions, Kalinowski and Pil\'sniak \cite{R. Kalinowski and M. Pilsniak} have defined the \textit{distinguishing index } $D'(G)$ of $G$ which is  the least integer $d$
such that $G$ has an edge coloring   with $d$ colors that is preserved only by a trivial
automorphism.  The distinguishing index and number of some examples of graphs was exhibited in \cite{Albert,R. Kalinowski and M. Pilsniak}. For 
 instance, $D(P_n) = D'(P_n)=2$ for every $n\geqslant 3$, and 
 $D(C_n) = D'(C_n)=3$ for $n =3,4,5$,  $D(C_n) = D'(C_n)=2$ for $n \geqslant 6$. It is easy to see that the value $|D(G)-D'(G)|$ can be large. For example $D'(K_{p,p})=2$ and $D(K_{p,p})=p+1$, for $p\geq 4$.  A symmetric tree, denoted by $T_{h,d}$, is a tree with a central vertex $v_0$, all leaves at
the same distance $h$ from $v_0$ and all the vertices which are not leaves with degree $d$. A bisymmetric tree, denoted by $T''_{h,d}$, is a tree with a central edge $e$, all leaves at the same distance $h$ from $e$ and all the vertices which are not leaves with degree $d$. The following theorem gives upper bounds for $D'(G)$ based on the maximum degree of $G$.  
 \begin{theorem}{\rm \cite{R. Kalinowski and M. Pilsniak, nord}}\label{11}
\begin{enumerate}
\item[(i)]  If $G$ is a connected graph of order $n\geq 3$, then
$D'(G) \leq \Delta(G)$, unless $G$ is $C_3$, $C_4$ or $C_5$.
\item[(ii)] Let $G$ be a connected graph that is neither a symmetric nor a bisymmetric
tree. If $\Delta(G)\geq 3$,  then
$D'(G) \leq \Delta(G)-1$, unless $G$ is $K_4$ or $K_{3,3}$.
\end{enumerate}
 \end{theorem}
 
 Also,  Kalinowski and Pil\'sniak  \cite{R. Kalinowski and M. Pilsniak} defined the \textit{distinguishing chromatic index}  $\chi'_D (G)$ of a graph $G$ as the least number $d$ such that $G$ has a proper edge labeling with $d$ labels that is preserved only by the identity automorphism
of $G$.
\begin{theorem}{\rm \cite{R. Kalinowski and M. Pilsniak}}\label{thm16} If $G$ is a connected graph of order $n \geq 3$, then
$\chi'_D (G) \leq \Delta (G) + 1$, except for four graphs of small order $C_4$, $K_4$, $C_6$, $K_{3,3}$. 
\end{theorem}
This theorem immediately implies the following interesting result. A proper edge labeling of $G$ with
$\chi'(G)$ colors is called minimal.
\begin{theorem}{\rm \cite{R. Kalinowski and M. Pilsniak}}\label{thm17} 
Every connected class 2 graph admits a minimal edge labeling that is not preserved by any nontrivial automorphism.
\end{theorem}
We need the following results. 

\begin{theorem}{\rm \cite{EJC}}\label{thm10}
If $T$ is a tree of order $n \geq 3$, then $D'(T) \leq \Delta(T)$. Moreover, equality is achieved if and only
if $T$ is either a symmetric  or a path of odd length. 
\end{theorem}
\begin{theorem}{\rm \cite{R. Kalinowski and M. Pilsniak}}\label{thm18}
 If $T$ is a tree of order $n \geq 3$, then
$\chi'_D(T ) = \Delta(T ) + 1$, if and only if $T$ is a bisymmetric tree.
\end{theorem}

In the next section, we compute the distinguishing chromatic index of certain graphs such as friendship and book graphs. More precisely, we present a table of results that shows the chromatic index,  the distinguishing index and the distinguishing chromatic index for  various families of
connected graphs. Also we  obtain a relationship between the chromatic distinguishing number of line graph $L(G)$ of graph $G$ and the chromatic distinguishing index of $G$. In Section 3, we study the distinguishing chromatic index of join and corona
product of two graphs.

\section{The  distinguishing chromatic index of certain graphs}
\begin{observation}
\begin{enumerate}
\item[(i)] For any graph $G$, $\chi'_D(G)\geq {\rm max}\{\chi'(G), D'(G)\}$.
\item[(ii)] If  $G$ has no non-trivial automorphisms, then $\chi'_D(G)= \chi'(G)$ and $D'(G)=1$. Hence $\chi'_D(G)$ can be much larger than $D'(G)$.
\end{enumerate}
\end{observation}

By this observation and Theorem \ref{thm16} we can conclude that for any connected graph $G$ of order $n\geq 3$ and maximum degree $\Delta$, $\chi'_D(G)$ is $\Delta$ or $\Delta+1$, except for $C_4$, $K_4$, $C_6$, $K_{3,3}$. In the latter case, $\chi'_D(G)=\Delta+2$. Hence, for any  connected graph $G$ we have $|\chi'_D(G) -\chi'(G)|\leq 2$, and equality is only achieved  for  $C_4$, $K_4$, $C_6$, and $K_{3,3}$.  By Theorem \ref{thm17}, it can be seen that $\chi'_D(G) =\chi'(G)$, for class 2 graphs. In the following theorem we present a family of class 1 graphs such that $\chi'_D(G) =\chi'(G)$.
\begin{theorem}\label{fixepoin}
Let $G$ be a class one graph, i.e., $\chi'(G)=\Delta(G)$. If there exists a vertex $x$ of $G$ for which $f(x)=x$ for all automorphisms $f$ of $G$, then $\chi'_D(G)=\Delta(G)$.
\end{theorem}
\proof
By  contradiction suppose  that $\chi'_D(G)=\Delta(G)+1$. Then, for any proper $\Delta(G)$-labeling $c$ of $G$, there exists a nonidentity automorphism $f$ of $G$ preserving the labeling $c$. By hypothesis, we have $f(x)=x$. Since the incident edges to $x$ have different labels, so $f$ fixes every adjacent vertex to $x$, because $f$ preserves the labeling $c$. By the same argument for every adjacent vertex to $x$, we can conclude that $f$ is the identity automorphism, which is a contradiction.\qed

By Theorems \ref{11}, \ref{thm16} and \ref{thm18}, we can characterize all connected graph with $\chi'_D(G)= D'(G)$.
\begin{theorem}
Let $G$ be a connected  graph of order $n\geq 3$ and maximum degree $\Delta$.
\begin{enumerate}
\item[(i)] There is  no  connected graph $G$ with $\chi'_D(G)= D'(G)= \Delta +2$.
\item[(ii)]  $\chi'_D(G)= D'(G)= \Delta +1$, if and only if $G\in \{ C_3,C_4,C_5\}$.
\item[(iii)]  $\chi'_D(G)= D'(G)= \Delta $, if and only if $G$ is a tree which is not a bisymmetric tree.
\end{enumerate}
\end{theorem}

\medskip
Here, we want to  obtain the distinguishing chromatic index of complete bipartite graphs. Before we obtain the distinguishing chromatic index of complete bipartite graphs we need the following information of \cite{fish}: A labeling with $c$ labels of the edges of a complete bipartite graph $K_{s,t}$ having parts $X$
 of size $s$ and $Y$ of size $t$ corresponds to a $t \times s$ matrix with entries from $\{ 1,\ldots , c\}$.
The $i, j$ entry of the matrix is $k$ whenever the edge between the $i$th vertex in $Y$ and
the $j$th vertex in $X$ has label $k$. We call this the bipartite adjacency matrix.  For edge labeled
complete bipartite graphs, the parts $X$ and $Y$ map to themselves if $|X|\neq |Y|$. In this
case, if $A$ is the bipartite adjacency matrix, then an automorphism corresponds to
selecting permutation matrices $P_Y$ and $P_X$ such that $A = P_Y AP_X$. If $|X|=|Y|$  then
we also have automorphisms of the form $A = P_Y A^T P_X$.   For any matrix with entries from $\{ 1,\ldots , c\}$ the degree of a column is a $c$-tuple
$(x_1,\ldots , x_c)$ with $x_i$ equal to the number of entries that are $i$ in the column.
\begin{theorem}{\rm \cite{fish}}\label{fact} 
Let $A$ be the adjacency matrix of a $c$-edge labeled complete bipartite graph.
\begin{enumerate}
\item[(i)] If there are two identical rows in $A$, then $A$ is not an identity labeling.
\item[(ii)] If $A$ is not square and if the columns of $A$ have distinct degrees and the rows
are distinct, then $A$ is an identity labeling. If $A$ is square, has distinct rows, distinct
column degrees and the multiset of column degrees is different from the multiset of
row degrees then $A$ is an identity labeling.
\end{enumerate}
\end{theorem}

\begin{theorem}\label{completebipartitechromaticindex}
The distinguishing chromatic index of complete bipartite graph $K_{s,t}$ where $s>t$, is $\chi'_D(K_{s,t})=s$.
\end{theorem}
\proof
Let $A$ be the following $s\times t$  adjacency matrix of a $s$-edge labeled complete bipartite graph,
\begin{equation*}
A=\left(\begin{array}{cccccc}
1&2&3&\cdots & s-1&s\\
s &1 &2 &\cdots &s-2 &s-1\\
s-1& s &1&\cdots &s-3&s-2\\
\vdots &\vdots &\vdots &\vdots &\vdots &\vdots \\
t&t+1&t+2&\cdots &t-2&t-1
\end{array}\right).
\end{equation*}
 Then it is clear that $A$ is a proper edge labeling. By Theorem \ref{fact} (ii), it can be concluded that $A$ is a distinguishing labeling. In fact, the rows of $A$ are distinct, and since the number of label $t+i-1 (mod ~s)$ in the $i$th column is one and in the $j$th column, $j\neq i$, is zero, so the columns have distinct degrees. Hence $\chi'_D(K_{s,t})=s$.\qed

 Before we prove the next result, we need the following preliminaries: 
  By the result obtained  by Fisher and Isaak \cite{fish} and independently by Imrich, Jerebic and Klav\v zar \cite{W Imrich} the distinguishing index of complete bipartite graphs is as  follows. 
\begin{theorem}{\rm \cite{fish, W Imrich}}\label{indcombipar}
	Let $p, q, r$ be integers such that $r \geq 2$ and $(r-1)^p <
	q \leq r^p$ . Then
	\begin{equation*}
	D'(K_{p,q}) =\left\{
	\begin{array}{ll}
	r & \text{if}~~ q \leq r^p - \lceil {\rm log}_r p\rceil - 1,\\
	r + 1 & \text{if}~~ q \geq r^p - \lceil {\rm log}_r p\rceil + 1.
	\end{array}\right.
	\end{equation*}
	If $q = r^p -\lceil {\rm log}_r p\rceil$ then the distinguishing index $D'(K_{p,q})$ is either $r$ or $r+1$ and can be computed recursively in $O({\rm log} ∗(q))$ time.
\end{theorem}

The friendship graph $F_n$ $(n\geqslant 2)$ can be constructed by joining $n$ copies of the cycle graph $C_3$ with a common vertex.
\begin{theorem}{\rm \cite{soltani2}}\label{indfriend}
Let $a_n=1+27n+3\sqrt{81n^2+6n}$. 
	For every $n\geq 2$, $$D'(F_n)=\lceil\frac{1}{3} (a_n)^{\frac{1}{3}}+\frac{1}{3(a_n)^{\frac{1}{3}}}+\frac{1}{3}\rceil.$$
\end{theorem} 

 The $n$-book graph $(n\geqslant 2)$  is defined as the Cartesian product $K_{1,n}\square P_2$. We call every $C_4$ in the book graph $B_n$, a page of $B_n$.   The distinguishing index of  Cartesian product of star $K_{1,n}$ with path $P_m$ for $m \geqslant≥ 2$ and $n \geqslant 2$ is $D'(K_{1,n}\square P_m) = \lceil \sqrt[2m-1]{n} \rceil$, unless $m=2$ and $n=r^3$ for some integer $r$. In the latter case $D'(K_{1,n}\square P_2) = \sqrt[3]{n}+1$, (\cite{gorz}).  Since $B_n=K_{1,n}\Box P_2$,  using this equality  we obtain the distinguishing index of book graph $B_n$.

\begin{theorem}
The entries in Table \ref{tabb} are correct.
\end{theorem}

\begin{table}[h]\label{tabb}
	\begin{center}
		\begin{tabular}{|c|c|c|c|c|}
			\hline
			& Graph $G$ & $\chi' (G)$& $D'(G)$ & $\chi'_D(G)$\\
			\hline 
			1. &  $K_{2n+1}$, $\Delta=2n$, $n=1,2$ & $2n+1$ & 3 & $2n+1$\\
			\hline 
			2. &  $ K_4$, $\Delta=3$ & 3& 3 & 5, Thm \ref{thm16}\\
			\hline 
			3. & $K_{2n}$, $\Delta=2n-1$, $n\geq 3$ &$2n-1$&  2 & $2n-1$\\
			\hline 
			4. &  $K_{2n+1}$, $\Delta=2n$, $n\geq 3$ & $2n+1$& 2 & $2n+1$, Thm \ref{thm17}\\\hline\hline
			5. &  $P_{2n}$, $\Delta=2$ & 2& 2 & 3\\\hline
			6. &  $P_{2n+1}$, $\Delta=2$ & 2& 2 & 2, Thm \ref{thm18}\\\hline\hline
			8. &  $C_4$, $\Delta=2$ & 2& 3 & 4, Thm \ref{thm16}\\ \hline
			9. &  $C_5$, $\Delta=2$ & 2& 3 & 3, Thm \ref{thm17}\\ \hline
			10. & $C_6$, $\Delta=2$ & 2& 2 & 4, Thm \ref{thm16}\\\hline
			11. & $c_{2n}$, $\Delta=2$, $n\geq 4$ & 2& 2 &3\\\hline
			12. &  $c_{2n+1}$, $\Delta=2$, $n\geq 3$  & 3& 2 & 3, Thm \ref{thm17}\\\hline\hline
			13. & $P$ Petersen, $\Delta=3$ & 4 & 3 & 4, Thm \ref{thm16} \\\hline
			14. &  bisymmetric tree, $\Delta$ & $\Delta$& $\Delta$, Thm \ref{thm10} & $\Delta+1$, Thm \ref{thm18}\\\hline
			15. &  tree $T\neq T''_{h,d}$,  $\Delta$ & $\Delta$& $\leq \Delta$, Thm \ref{thm10} & $\Delta$, Thm \ref{thm18}\\\hline
			16. &  $K_{3,3}$, $\Delta=3$ & $3$& 3 &5, Thm \ref{thm16}\\\hline
			17. &  $K_{n,n}$, $\Delta=n$, $n\geq 4$ & $n$& 2 & $n+1$\\\hline
			18. &  $K_{n,m}$, $m > n$, $\Delta=m$ & $m$& Thm \ref{indcombipar} & $m$\\\hline
			19. &  $F_{n}$, $n\geq 3$, $\Delta=2n$ & $2n$& Thm \ref{indfriend} & $2n$\\\hline
			20. &  $B_n$, $n\geq 2$, $\Delta=n+1$ & $n+1$& Thm in \cite{gorz} & $n+1$\\\hline
		\end{tabular}
		\caption{Tabel of results for $\chi$, $D'$ and $\chi'_D$.}
	\end{center}
\end{table}

\proof The chromatic index and the distinguishing index for the classes of graphs given in this table are well-known, we justify the entries in the last column.

\medskip
\noindent\textbf{Paths of even order.}
The labeling that uses label 2 for one end-edge and label $1$ for the remaining edges is distinguishing, however, it is not a proper labeling and any proper labeling using two labels is not distinguishing, so $\chi'_D(P_{2n})\geq 3$. A $3$-labeling that is proper and distinguishing is achieved by using $1$ for an end-edge and alternating 2's and 3's for the remaining edges, thus $\chi'_D(P_{2n})\leq 3$. Thus the result follows.

\medskip
\noindent\textbf{Cycles of even order $2n$, where $n\geq 4$.} Let the consecutive edges of $C_{2n}$ be $e_1,e_2,\ldots , e_{2n}$. Using label 3 for edges $e_1$  and $e_4$, label 2 for edges $e_i$ where $i\neq 1$ is odd and label 1 for edges $e_i$ where $i\neq 4$ is even, we get $\chi'_D(C_{2n}) \leq 3$ for $n\geq 4$. All proper 2-labelings of edges of $C_{2n}$ have label preserving automorphisms, thus  $\chi'_D(C_{2n})> 2$ for all $n$. Therefore  $\chi'_D(C_{2n}) = 3$, where $n\geq 4$.

\medskip
\noindent\textbf{Friendship and book graphs.} It is clear that each proper $2n$-labeling of edges of $F_n$ is distinguishing and so $\chi'_D(F_n) = 2n$, by Theorem \ref{fixepoin}.  For the book graph $B_n$, we  present a proper $(n+1)$-distinguishing edge labeling. Let $v_0$ and $w_0$ be two vertices of degree $n+1$ of $B_n$, and $v_1, \ldots , v_n$ be the adjacent vertices to $v_0$, and $w_1,\ldots , w_n$ be the adjacent vertices to $w_0$, such that $v_0v_iw_iw_0$ for $1\leq i \leq n$ are pages of $B_n$. We label  the edges $v_0v_i$, $v_iw_i$ and $w_iw_0$ with labels $i$, $i+2$ and $i+1$ mod $n$, respectively, for any $i$, $1\leq i \leq n$. Also, we label the edge $v_0w_0$ with label $n+1$. It can be seen that our labeling is a proper $(n+1)$-distinguishing labeling.

\medskip
\noindent\textbf{Complete graphs of even order.} It is known that we can partition the edge set of $K_{2n}$ to $2n-1$ sets, each set contains an $n$-element perfect matching, say $M_1, \ldots , M_{2n-1}$. If we label the edges of $n$-element perfect matching $M_i$ with label $i$, for any $1\leq i \leq 2n-1$, then it can be seen that this labeling is proper. We claim that this labeling is distinguishing. If $f$ is an automorphism of $K_{2n}$ preserving the labeling, then $f$ fixes the set $M_i$, for any $1\leq i \leq 2n-1$, setwise. If $f$ is a nonidentity automorphism, then without loss of generality we can assume that $f(v_1)=v_2$. Since $f$ preserves the labeling so $f(v_2)= v_1$. We can suppose that there exists a vertex $x$ of $K_{2n}$ such that the edges $\{v_1,x\}$ and $\{v_2,f(x)\}$ are not in the same perfect matching. Thus the labels of edges $\{v_1,x\}$ and $\{v_2,f(x)\}$ is different, while $f$ maps these two edges to each other, which is a contradiction. Then, the identity automorphism is the only automorphism of $K_{2n}$ preserving the labeling, and hence the proper edge labeling is  distinguishing, and so $\chi'_D(K_{2n})=2n-1$.

\medskip
\noindent\textbf{Complete bipartite graph $K_{n,n}$, $n\geq 4$.} We can partition the edge set of $K_{n,n}$ to $n$ perfect matching $M_1,\ldots , M_n$, each of $M_i$ contains $n$ edges. For every value of $n$, the labeling that uses label $i$ for all edges in $M_i$, $1\leq i\leq n$, is a proper labeling. Also, every proper labeling of $K_{n,n}$ partitions the edges of $K_{n,n}$ to $n$ sets $N_1, \ldots , N_n$ such that  each of $N_i$ is a perfect matching of $K_{n,n}$ and all edges in $N_i$ have the same label, and different from the label of  edges in $N_j$ for every $i,j\in \{1, \ldots , n\}$ where $i\neq j$. However, for every proper labeling of  $K_{n,n}$ we can find a nonidentity automorphism of  $K_{n,n}$ preserving the labeling, thus $\chi'_D(K_{n,n}) > n$. Now the result follows from Theorem \ref{thm16}.\qed

\medskip
In sequel, we want to obtain a relationship between  the chromatic distinguishing number  of line graph $L(G)$ of graph $G$ and  the chromatic distinguishing index  of $G$. For this purpose, we need more information about automorphism group of $L(G)$. For a simple graph $G$, we recall that  the \textit{line graph}  $L(G)$ is a  graph whose vertices are edges of $G$ and where two edges $e, e' \in   V (L(G)) = E(G)$ are adjacent if they share an endpoint in common.  Let  $\gamma_G : {\rm Aut} (G) \rightarrow {\rm Aut} (L(G))$ be given by $(\gamma_G \phi)(\{u, v\}) = \{\phi(u), \phi(v)\}$ for every $\{u, v\} \in  E(G)$. In \cite{Sabidussi}, Sabidussi proved the following theorem which we will use
throughout.
\begin{theorem}{\rm \cite{Sabidussi}}\label{autlinegraph}
 Suppose that $G$ is a connected graph that is not $P_2, Q$, or $L(Q)$ (see Figure \ref{autline}). Then 
$G$ is a group isomorphism, and so ${\rm Aut}(G) \cong {\rm Aut}(L(G))$.
\begin{figure}
	\begin{center}
		\includegraphics[width=0.55\textwidth]{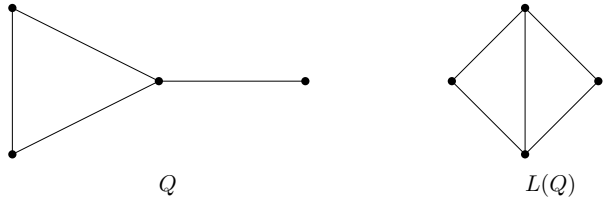}
		\caption{\label{autline} graphs $Q$ and $L(Q)$ of Theorem \ref{autline}.}
	\end{center}
\end{figure}
\end{theorem}

%Now we are ready to obtain a relationship between  the chromatic distinguishing number  of line graph $L(G)$ of graph $G$ and  the chromatic distinguishing index  of $G$.

\begin{theorem}
Suppose that $G$ is a connected graph of order $n\geq 3$ that is not $Q$ and $L(Q)$. Then $\chi_D(L(G))= \chi'_D(G)$.
\end{theorem}
\proof   First, we show that $\chi_D(L(G))\leq  \chi'_D(G)$. For this purpose, let $c:E(G)\rightarrow \{1, \ldots , \chi'_D(G)\}$ be an edge proper distinguishing labeling of $G$. We define $c': V(L(G)) \rightarrow \{1, \ldots , \chi'_D(G)\}$ such that $c'(e) = c(e)$ where $e\in V(L(G))  = E(G)$. By the following steps we show that the vertex labeling $c'$ is  proper distinguishing labeling. 
\begin{enumerate}
\item[Step 1)] The vertex labeling $c'$ is  proper. If $e$ and $e'$ are two adjacent vertices of $L(G)$ with $c'(e)=c'(e')$, then it means that $e$ and $e'$  are incident edges of $G$ with $c(e)=c(e')$, which is impossible. Thus $c'$ is  a proper labeling.
\item[Step 2)] The vertex labeling $c'$ is a distinguishing vertex labeling of $L(G)$, because if $f$ is an automorphism of $L(G)$ preserving the labeling, then $c'(f(e))= c'(e)$, and hence $c(f(e))= c(e)$ for any $e\in E(G)$. On the other hand $f = \gamma_G \phi$ for some automorphism $\phi$ of $G$, by Theorem \ref{autlinegraph}. Thus from  $c(f(e))= c(e)$ for any $e\in E(G)$, we can conclude that $c(\gamma_G \phi (e))= c(e)$ and so $c(\{\phi(u), \phi(v)\})= c(\{u,v\})$ for every $\{u,v\}\in E(G)$. This means that $\phi$ is an automorphism of $G$ preserving the labeling $c$, so $\phi$ is the identity automorphism of $G$. Therefore $f$ is the identity automorphism of $L(G)$, and hence $\chi_D(L(G))\leq  \chi'_D(G)$.  
\end{enumerate}

By a similar argument we can prove that $\chi_D(L(G))\geq \chi'_D(G)$, and so the result follows.\qed

\begin{figure}
	\begin{center}
		\includegraphics[width=0.4\textwidth]{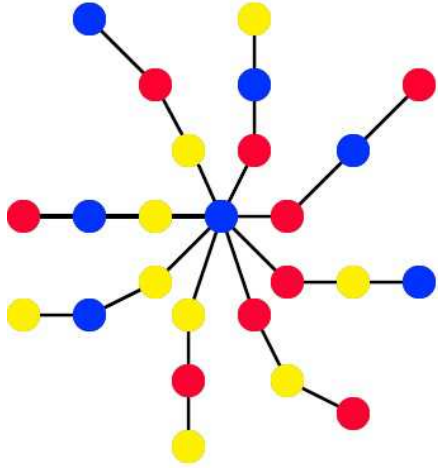}
		\caption{A $3$-proper distinguishing vertex labeling of $G_8$.}\label{fig1}
	\end{center}
\end{figure}

\medskip
Now we state a difference between the distinguishing labeling and chromatic distinguishing labeling of graphs.
\begin{remark}
Despite the distinguishing labeling that $D'(G) \leq D(G) +1$ for all connected graph $G$, it may be happen that $\chi'_D(G) > \chi_D(G) +1$. For instance, let $G_n$ be the graph obtained from $K_{1,n}$ by replacing each edge  with a path of length three. It can be computed that $ \chi'_D(G_8)=8$, while $\chi_D(G_8)=3$, see figure \ref{fig1}.  
\end{remark}

We end this section by proposing  the following problem.
\begin{problem}
Characterize all connected graphs  with $\chi'_D(G)=\chi'(G)=\Delta(G)$.
\end{problem}

\section{Results for join and corona products}

In this section we study the distinguishing chromatic index of join and corona product of graphs. We start with join of graphs.  The graph $G = (V, E)$ is the \textit{join} of two graphs $G_1 = (V_1, E_1)$ and
$G_2 = (V_2, E_2 )$, if $V = V_1 \cup V_2$ and $E = E_1 \cup E_2 \cup \{uv | u \in V_1, v \in V_2\}$
and denoted by $G = G_1 + G_2$.

\begin{theorem}\label{join}
If   $G$ and $H$ are two connected graphs of orders $n,m\geq 3$,  respectively, then
\begin{equation*}
{\rm max}\{\Delta(H)+n, \Delta(G)+m\}\leq \chi'_D(G+H)\leq {\rm max}\{\chi'_D(G),\chi'_D(H)\}+\chi'_D(K_{n,m}).
\end{equation*}
\end{theorem}
\proof  To prove the left inequality, it is sufficient to know that $\Delta(G+H)= {\rm max}\{\Delta(H)+n, \Delta(G)+m\}$ and  $\chi'_D(G+H) \geq \Delta(G+H)$. For the right inequality, we first set $M:=  {\rm max}\{\chi'_D(G), \chi'_D(H)\}$. We label the edge set of graph $G$ (resp. $H$) with labels $\{1,\ldots ,\chi'_D(G) \}$ (resp. $\{1,\ldots ,\chi'_D(H)\}$)  in a proper distinguishing way.  If $V(G)=\{v_1,\ldots , v_n\}$ and $V(H)=\{w_1,\ldots , w_m\}$, then we label the middle edges $\{v_iw_j~:~1\leq i\leq n, 1\leq j \leq m\}$, exactly the same as a proper distinguishing labeling of the complete bipartite graph $K_{|V(G)|, |V(H)|}$ with labels $\{M+1, \ldots , M+ \chi'_D(K_{|V(G)|, |V(H)|})\}$. Since the graphs $G$ and $H$ have a proper labeling, so this labeling of edges of $G+H$ is proper, regarding to the label of middle edges. To show this labeling is distinguishing, we suppose that $f$ is an automorphism of $G+H$ preserving the labeling. Then, with respect to the label of middle edges, it can be concluded that the restriction of $f$ to the vertices of $G$ (resp. $H$) is an automorphism of $G$ (resp. $H$) preserving the labeling. Since the graphs $G$ and $H$ have been labeled  distinguishingly, so the restriction of $f$ to vertices of $G$ and $H$ is the identity automorphism of $G$ and $H$, respectively. Therefore, this labeling is distinguishing. \qed

The lower and upper bound of Theorem \ref{join} are sharp. For example, we consider $P_{2n+1}+P_{2m+1}$ where $m>n$, with $\chi'_D(P_{2n+1}+P_{2m+1})= 2m+3$. In fact, we can label the edge set of $P_{2n+1}$ and $P_{2m+1}$ with labels 1 and 2 in a proper distinguishing way, and label the remaining incident edges to each vertex of $P_{2n+1}$ with distinct labels $3,4,\ldots , 2m+3$ such that  the  incident edges to each vertex of $P_{2m+1}$ in $P_{2n+1}+P_{2m+1}$ have distinct labels. 

\medskip
Now we want to obtain a lower and upper bound for the distinguishing chromatic index of corona product. 
The \textit{corona product}  $G\circ H$ of two graphs $G$ and $H$ is defined as the graph obtained by taking one copy of $G$ and $\vert V(G)\vert $ copies of $H$ and joining the $i$-th vertex of $G$ to every vertex in the $i$-th copy of $H$.  If  $ G$ and $H$ be two connected  graphs such that $G\neq K_1$, then  there is no vertex in the copies of $H$ which has the same degree as a vertex in $G$. Because if  there exists a vertex $w$ in one of the copies of $H$  and a vertex $v$ in $G$ such that ${\rm deg}_{G\circ H} (v)={\rm deg}_{G\circ H} (w)$, then  ${\rm deg}_G(v)+|V(H)|={\rm deg}_H(w)+1$. So we have ${\rm deg}_H(w)+1 > |V(H)|$, which is a contradiction. Hence, we can state the following lemma.

\begin{lemma}\label{lem2}
Let   $ G$ and $H$ be two connected  graphs such that $G\neq K_1$. If  $f$ is an arbitrary automorphism  of $G\circ H$, then the restriction of $f$ to the vertices of copy $G$ (resp. $H$) is an automorphism of $G$ (resp. $H$).
\end{lemma}

\begin{theorem}\label{corona}
If  $G$ and $H$ are two connected graphs of orders $n,m\geq 3$, respectively,  then  
\begin{equation*}
{\rm max}\{\Delta(G)+m, \Delta(H)+1\}\leq \chi'_D(G\circ H)\leq {\rm max}\{\chi'_D(G), \chi'_D(H)\} +m.
\end{equation*}
\end{theorem}
\proof  The proof of the left inequality is exactly the same as Theorem \ref{join}. To prove the right inequality, we first set $M:=  {\rm max}\{\chi'_D(G), \chi'_D(H)\}$. We label the edge set of graph $G$ (resp. $H$) with labels $\{1,\ldots ,\chi'_D(G) \}$ (resp. $\{1,\ldots ,\chi'_D(H)\}$)  in a proper distinguishing way.  If $V(G)=\{v_1,\ldots , v_n\}$ and $V(H)=\{w_1,\ldots , w_m\}$, then we label the middle edge $v_iw_j$ with label $M+j$ for any $1\leq i\leq n$ and $1\leq j \leq m$. Since the graphs $G$ and $H$ have been labeled  properly, so this labeling of edges of $G\circ H$ is proper, due to the label of middle edges. To show this labeling is distinguishing, we suppose that $f$ is an automorphism of $G\circ H$ preserving the labeling. Then, by Lemma \ref{lem2}, it can be concluded that the restriction of $f$ to the vertices of $G$ (resp. $H$) is an automorphism of $G$ (resp. $H$) preserving the labeling. Since the graphs $G$ and $H$ have been labeled  distinguishingly, so the restriction of $f$ to vertices of $G$ and $H$ is the identity automorphism of $G$ and $H$ respectively. Therefore, this labeling is distinguishing. \qed

The lower and upper bound of Theorem \ref{corona} are sharp. For instance, we consider $K_{2n}\circ K_{2m+1}$ where $3\leq m < n$. It can be easily computed that  $\chi'_D(K_{2n}\circ K_{2m+1})= 2m+2n$.  Now, since ${\rm max}\{\Delta(K_{2n})+2m+1, \Delta(K_{2m+1})+1\}=  2m+2n$ and ${\rm max}\{\chi'_D(K_{2n}), \chi'_D(K_{2m+1})\} +2m+1= 2m+2n$, so the bounds of Theorem \ref{corona} are sharp.

\medskip
We end this paper  by a remark on the distinguishing chromatic index of join and corona product of graphs $G$ and $H$ where $G=K_1$ or $K_2$.
\begin{remark}
\begin{enumerate}
\item[(i)] If $G$ and $H$ are two connected graphs such that $G=K_1$ or $K_2$, then the maximum degree of $G+H$ is $|V(H)|$ or $|V(H)|+1$, respectively, and so $\chi'_D(G+H) \leq |V(H)|+2$, except for $G=H= K_2$. In the latter case, $\chi'_D(K_2+K_2) =5$.
\item[(ii)] It is clear that $K_1 \circ H = K_1 +H$, and so $\chi'_D(K_1 \circ H) = \chi'_D(K_1 +H)$.  If $G=K_2$, then $\Delta (K_2 \circ H)= |V(H)|+1$, and so $\chi'_D(K_2 \circ H) \leq |V(H)|+2$.
\end{enumerate}
\end{remark}
 
\end{document}